\documentclass[12pt]{article}%
\usepackage{graphicx}
\usepackage{algorithm2e}
\usepackage{amsmath, amsfonts, amsthm, amssymb}
\usepackage{geometry}
\usepackage{subcaption}
\usepackage{caption}%
\usepackage{amsmath}%
\usepackage{amsfonts}%
\usepackage{amssymb}
\let\epsilon\varepsilon
\newcommand\hide[1]{}
\makeatother
\theoremstyle{definition}

\SetKwComment{Comment}{/* }{ */}
\SetKwFor{For}{for (}{) $\lbrace$}{$\rbrace$}
\RestyleAlgo{ruled}
\begin{document}

\title{DeepPropNet - A Recursive Deep Propagator Neural Network for Learning
Evolution PDE Operators}
\author{Lizuo Liu, Wei Cai \thanks{Department of Mathematics, Southern
Methodist University, Dallas, TX 75252. email: cai@smu.edu. Date: February 27, 2022.}}
\maketitle

\begin{abstract}
In this paper, we propose a deep neural network approximation to the
evolution operator for time dependent PDE systems over long time period by
recursively using one single neural network propagator, in the form of POD-DeepONet with built-in causality feature,
for a small time interval. The trained DeepPropNet of moderate size is
shown to give accurate prediction of wave solutions over the whole time interval.

\end{abstract}

\section{The Problem of Learning PDE Evolution Operator}

\label{sec:evoop}

Consider a second order evolution system%
\begin{align}
u_{tt} &  =\mathcal{L}u+f(x,t),\label{evol}\\
u(x,0) &  =u_{0}(x),\quad x\in R^{1},\nonumber\\
u_{t}(x,0) &  =v_{0}(x).
\end{align}
%
where $u$ could be a scalar or vector. For a model problem, we will consider
the inhomogeneous scalar wave (d'Alembert) equation
\begin{equation}
\frac{\partial^{2}u}{\partial t^{2}}-c^{2}\left(  x,t\right)  \frac
{\partial^{2}u}{\partial x^{2}}=f(x,t),\quad x\in\mathbb{R},\quad t\geq0,
\label{1dwave}
\end{equation}
with the source term $f$ that is
compactly supported on a bounded space-time domain $Q=S\times\left[
0,T \right]  $ where $S\subset\mathbb{R}$. This means that the source term
$f(x,t)$ differs from zero only on $S$ and operates for the
limited time interval $\left[  0,T \right]  $. When $c\left(  x,t\right)
=\text{const, t}$he solution to the Cauchy problem $\left(  \ref{1dwave}%
\right)  $ is given by:%

\begin{equation}
u(x,t)=\frac{1}{2}(u_{0}(x-ct)+u_{0}(x+ct))+\frac{1}{2c}\int_{x-ct}%
^{x+ct}v_{0}(\xi)d\xi+\frac{1}{2c}\int_{0}^{t}\int_{x-c(t-\tau)}^{x+c(t+\tau
)}f(\xi,\tau)d\xi d\tau.\label{Dalemb}%
\end{equation}


The solution $u(x,t)$ can be viewed through an evolution operator $P(u_{0},v_{0,}f(x,s),0\leq s\leq
t)$, which maps the initial conditions and the source term into the solution.  Learning such a map between functions has
been actively studied recently with various types of operator learning methods, including DeepONet \cite{cdeepOnet} and Fourier Neural Operator \cite{fno}. The focus of this paper is to find an efficient way
to learn this evolution operator with moderate size neural network for large time $t$.
For lack of a precise term, borrowing the term from quantum
mechanics for the Green's function propagator \cite{shankar}, we shall name the
operator $P$ as the propagator for the evolution system. From (\ref{Dalemb}),
it is clear that if we like to train a neural network operator for large time
t, the size of the network will grow for increasing time t. And, the amount of
information to be input into a network will increase dramatically as t grows
as well. We will propose a recursive propagator formulation for
the evolution operator network.

First, the solution time interval $[0,T]$ will be divided into $N$ smaller
subintervals%
\begin{equation}
t_{0}=0<t_{1}<\cdots<t_{i}<\cdots t_{N},t_{i}=i\Delta t,\quad\Delta
t=T/N,\label{mesh}%
\end{equation}
and for $t_{i}\leq t\leq t_{i+1},$ the solution is given by the propagator
with initial condition of the solution $u_{i}$ and its velocity $v_{i}%
=\overset{.}{u}_{i}$, i.e.,%
\begin{equation}
u(x,t)=P(u_{i},v_{i},f(x,s),t_{i}\leq s\leq t),\label{propag2}%
\end{equation}
where the initial condition $u_{i}$ would have been given by the propagator
for the time block $t_{i-1}\leq t\leq t_{i}$.

The propagators in (\ref{propag2}) will be approximated by a single neural
network in the form of DeepONet structure \cite{deeponet}  with modification for time
causality, namely,%
\begin{equation}
P(u_{i},v_{i},f(x,s),t_{i}\leq s\leq t)\sim P_{\theta}(u_{i},v_{i}%
,f(x,s),t_{i}\leq s\leq t),0\leq i\leq N-1.\label{propag3}%
\end{equation}

Therefore, the propagator $P_{\theta}(u_{0}(x),v_{0}(x),f(x,t))$ will
be trained to map the initial and force function data into the solution
$u(x,t),t_{0}\leq t\leq t_{1}$. Moreover, this same propagator $P_{\theta}$
will be trained to approximate the solution for time periodic $t_{1}\leq t\leq
t_{2}$ where the initial condition at $t_{1}$ can be computed with the
propagator for the previous time interval $[t_{0},t_{1}]$. This procedure will
be used recursively until we have trained the same propagator for the last
time interval $[t_{N-1},t_{N}]$. Taken all together, we arrive at a propagator
neural network for the whole time interval $[0,T]$ where the building block is
the single propagator of moderate size $P_{\theta}(u_{0}(x),v_{0}(x),f(x,t)),t_{0}\leq t\leq t_{1})$. By controlling the size of the
$\Delta t$, the size of this propagator can be easily controlled for
efficiency as well as accuracy. As the evolution PDE system has to observe the
causality of the physical system, the DeepONet framework will be modified to
include the causality, a previously proposed causality DeepONet in the study
of dynamics system for modeling building response to seismic waves
\cite{cdeepOnet} will be used for this purpose.

\medskip

The rest of the paper is organized as follows. In section \ref{sec:cpod}, we will review
the DeepONet \cite{cdeepOnet} with time causality and extension with proper orthogonal decomposition (POD)
approach for efficient treatment of spatial dependence of the solution. Section \ref{sec:dp} will give the
algorthmic outline of the DeepPropNet and numerical results of the DeepPropNet
will be presented in Section \ref{sec:results}. Finally, Section \ref{sec:future} gives a conclusion and some
future work.

\section{DeepONet with time causality and spatial POD}
\label{sec:cpod}
A causality DeepONet was proposed in \cite{cdeepOnet} to handle the time
casuality in dynamics system and was shown to be very effective to predict the
seismic response of building. Here, we will just present the final form and
for details please refer to \cite{cdeepOnet}.

\textbf{Causality-DeepONet:} A DNN representation of an operator
$\mathcal{G}(f)(t)$ for any continuous function $f(t)$ with retarded response
for $t\in K_{2}=[0,T]\subset$ $\mathbb{R}$ is given as
\begin{equation}
\mathcal{G}(f)(t) \sim\sum_{k=1}^{N}\sum_{i=1}^{M}c_{i}^{k}\sigma_{b}\left(
\sum_{j=1}^{\left\lceil \frac{t}{h}\right\rceil }\xi_{i,m-\left\lceil \frac
{t}{h}\right\rceil +j}^{k}f\left( s_{j}\right)  \right)  \cdot\sigma
_{trk}\left(  \omega_{k}\cdot t+\zeta_{k}\right),  \label{operator_appro}%
\end{equation}
where $\{s_{j}\}_{j=1}^{m}\subset K_{1}=[0,T]\subset\mathcal{X},$coefficents
$c_{i}^{k},\xi_{ij}^{k},\omega_{k},\zeta_{k}$-all independent of continuous
functions $u$ $\in V\subset C(K_{1})$ and $t$.


To handle the spatial dependence of solution $u(x,t)$ for the evolution
system, we will adopt the idea of the
POD-DeepONet\cite{ModelReductionNeural2021a,luComprehensiveFairComparison2021}
which assumes there is a set of global basis for the targeted output which
could be found by SVD, and the trunk net will be replaced by these basis, then the neural
network is learning the mapping between values at sensors to
the singular values if we consider using SVD as an example. Thus it may lose
the ability to approximate arbitrary values of the operator in the computational domain
given the values at sensors. The philosophy of POD-DeepONet is to find a set of basis globally,
thus compress the memory and computation needed, which is the crux for
learning the operator of the high dimensional problem$\left(  d\gg O\left(  1
\right)  \right) $.

Lu et al.\cite{luComprehensiveFairComparison2021} discussed the POD-DeepONets
and proposed a time marching scheme for Fourier Neural Operator specifically
to reduce the dimension of output, and they also proposed the modified
DeepONet with feature expansion by feeding historical states of signal to
trunk net as features. Bhattacharya et al.\cite{ModelReductionNeural2021a}
proposed that they reduce the dimension of input and output by PCA and then
learning the mapping between the reduced space by a neural network, and they
also discussed the approximation quality needed to get good mappings. Meuris
et al.\cite{MachinelearningCustommadeBasis2021} presented a procedure that
learning the data-driven basis functions harnessing the DeepONet machinery and
the learned basis function will be used by classical methods as custom basis
to achieving high accurary for computation of arbitrary complex domain.

As for a problem with causality, The POD-DeepONet's philosophy guides us
either do the SVD for the whole outputs regardless the difference of temporal
variables and spatial variables, or do specific SVD timestep by timestep to
keep the causality for which we need to keep the basis for all time steps, or
consider solution at each time step for each case as an independent target and
find a common basis w.r.t spatial dimensions for all time steps \& all cases. These ideas either burns
the high memory cost, or destroy the causality in a brutal-force way.
Therefore, a modification of POD-DeepONet for problem with causality is
crucial. We follow the idea of POD-DeepONet but only construct the basis of
spatial domain explicitly and utilize the causality DeepOnet to handle the temporal-dependent coefficients of each spatial basis.

In the following, we consider the 1d wave equation (\ref{1dwave} and assume the right hand side
has form
\begin{equation}
f\left(  x,t \right)  = a_{0}\left(  t \right)  + \sum_{n=1}^{N} a_{n}\left(
t \right)  \cos\left(  2 n \pi x \right)  + b_{n}\left(  t \right) \sin\left(
2 n \pi x \right) ,\label{1drhs}%
\end{equation}
and the solution has a similar form
\begin{equation}
u\left(  x,t \right)  = \psi_{0}\left(  t \right)  + \sum_{m=1}^{M} \psi
_{m}\left(  t \right)  \cos\left(  2 m \pi x \right)  + \phi_{m}\left(  t
\right) \sin\left(  2 m \pi x \right) ,\label{1dsol}%
\end{equation}

Next, we modify the Causality-DeepONet $\left(  \ref{operator_appro}
\right) $ by
\begin{equation}
\mathcal{G}(\overrightarrow{f})(x,t) \sim\overrightarrow{\sigma}_{br}\left(
\overrightarrow{f} \right) \odot\overrightarrow{\sigma}_{trk}\left(  t
\right)  \cdot\overrightarrow{\sigma}_{basis}\left(  x \right),  \label{cpod}%
\end{equation}
where $\odot$ is the elementwise multiplication, $\cdot$ is the inner
product,
\[
\overrightarrow{f}\left(  s \right)  = \left[  a_{0}\left(  s \right)
,a_{1}\left(  s \right) ,b_{1}\left(  s \right) ,\ldots,a_{N}\left(  s \right)
,b_{N}\left(  s \right)  \right] ,
\]
\begin{align}%
\begin{split}
\overrightarrow{\sigma}_{br, i}\left(  \overrightarrow{f} \right)  ={} &
\sigma_{b}\left(  \sum_{j=1}^{\left\lceil \frac{t}{h}\right\rceil }%
\xi_{i,m-\left\lceil \frac{t}{h}\right\rceil +j}^{0}a_{0}\left(  s_{j}\right)
\right. \\
& \left.  + \sum_{n=1}^{N} \sum_{j=1}^{\left\lceil \frac{t}{h}\right\rceil
}\left(  \xi_{i,m-\left\lceil \frac{t}{h}\right\rceil +j}^{n,a}a_{n}\left(
s_{j}\right)  + \xi_{i,m-\left\lceil \frac{t}{h}\right\rceil +j}^{n,b}%
b_{n}\left(  s_{j}\right)  \right) \right) ,\label{branch_cpod}%
\end{split}
\\%
\begin{split}
\overrightarrow{\sigma}_{trk,i}\left(  t \right)  ={} &  \sigma_{trk}\left(
\omega_{i}\cdot t+\zeta_{i}\right) ,\label{trunk_cpod}%
\end{split}
\\%
\begin{split}
\overrightarrow{\sigma}_{basis}\left(  x \right)  ={} &  \{ 1, \cos(2\pi x),
\sin(2\pi x), \cdots, \cos(2\pi N x), \sin(2\pi N x) \}. \label{basis}%
\end{split}
\end{align}


The activation functions of $\sigma_{b}$ and $\sigma_{trk}$ are ReLU. Note the
modes number $M$ in the CPOD-DeepONet(Causality POD-DeepONet) \(\left( \ref{1dsol} \right)\) need to be
greater than or equal to the number of given modes \(N\) of right hand side in equation \(\left( \ref{1drhs}  \right)\). The algorithm of the CPOD-DeepONet is shown in algorithm
\ref{alg:CPOD-DeepONet}. To contain the memory efficiency and following the
fact that equation $\left(  \ref{1dsol} \right) $ has separation of variables,
we will do outer product rather than computing point by points as shown in the
last line of algorithm. This is one of the crux to handle higher dimensional problem.

\begin{algorithm}
\caption{The algorithm of Causality POD-DeepONet}\label{alg:CPOD-DeepONet}
\SetKwInOut{Input}{Input}
\SetKwInOut{Output}{Output}
\Input{$t$, \({x}\), \(\left[ a_0\left( s \right), a_1\left( s \right), b_1\left( s \right),\ldots ,a_{N}\left( s \right),b_{N}\left( s \right) \right] \)}
\Output{$u\left( {x},t \right)$}
\(y_t \gets \text{CDeepONet}\left( t, \left[ a_0\left( s \right), a_1\left( s \right), b_1\left( s \right),\ldots ,a_{N}\left( s \right),b_{N}\left( s \right) \right] \right)\) \Comment*[r]{CDeepONet is the causality DeepONet. The output \(y_t\) has size of \([\text{batchsize} \times  \text{Nt} \times \text{NBx}]\)}
$y_x \gets \left[ 1, \cos\left( 2\pi {x} \right), \sin\left( 2\pi {x} \right), \ldots \cos\left( 2 M \pi {x} \right), \sin\left( 2 M {x} \right)\right] $\Comment*[r]{The output \(y_x\) has size of \([\text{batchsize} \times  \text{Nx} \times \text{NBx}]\)}
$u\left( x,t \right) \gets \text{Einsum}\left( \text{'bxn,btn} \to \text{ bxt'}, y_{x}, y_{t} \right) $\Comment*[r]{The Einsum is the Einstein summation convention, \(\left( \text{'bxn,btn} \to \text{ bxt'} \right)\) means the index change in the Einstein summation convention. }
\end{algorithm}

\textbf{Loss function} Given batch size $\mathcal{N}$ for training process and
the total number of the test records $N$, the loss function is defined as
\begin{equation}
\mathcal{L}oss\left(  \theta\right)  =\frac{1}{\mathcal{N}}\sum_{i=1}%
^{\mathcal{N}}\sum_{j=1}^{N_{t}}\sum_{k=1}^{N_{{x}}} \left(
\mathcal{G}\left(  \overrightarrow{f} \right) \left(  {x}_{k}%
,t_{j}\right)  -y_{ijk}\right)  ^{2}, \label{loss}%
\end{equation}
where $N_{t}$ is the number of time step and $N_{x}$ is the number of points
on ${x}$ direction.

To evaluate the training process, the mean of the relative L2 error is
considered
\begin{equation}
\mathcal{L}oss_{2}\left(  \theta\right)  = \frac{1}{\mathcal{N}}\sum
_{i=1}^{\mathcal{N}}\frac{ \left\|  \mathcal{G}\left(  \overrightarrow{f_{i}}
\right) \left(  {x},t \right)  - u_{i}\left(  {x},t
\right)  \right\| _{2}}{\left\|  u_{i}\left(  {x},t \right)
\right\| _{2} }.
\end{equation}
The relative L2 error in a complete epoch is defined as
\begin{equation}
\mathcal{L}^{\mathcal{R}}_{\text{train}} = \frac{1}{B} \sum_{k=1}^{B} \frac
{1}{\mathcal{N}}\sum_{i=1}^{\mathcal{N}}\frac{ \left\|  \mathcal{G}\left(
\overrightarrow{f_{i}} \right) \left(  {x},t, \theta^{\left(  k
\right) }\right)  - u_{ik}\left(  {x},t \right)  \right\| _{2}%
}{\left\|  u_{i}\left(  {x},t \right)  \right\| _{2} },
\label{L_2_train}%
\end{equation}
where $B$ is the number of batches, $\theta^{\left(  k \right)  }$ means the
parameters of neural network at $k$-th batch.

Similarly, we define the relative L2 error for the testing dataset
\begin{equation}
\mathcal{L}_{\text{test}}^{\mathcal{R}}=\frac{1}{{N}}\sum_{i=1}^{N} \frac{
\left\|  \mathcal{G}\left(  \overrightarrow{f_{i}} \right) \left(
{x},t \right)  - u_{i}\left(  {x},t \right)  \right\|
_{2}}{\left\|  u_{i}\left(  {x},t \right)  \right\| _{2} }.%
\end{equation}
Note the $N$ is the total number of test cases.

\bigskip

\noindent$\bullet$ \textbf{Numerical Performance of Causality DeepONet with
POD}

\medskip

\begin{figure}[ptbh]
\makebox[\textwidth][c]{\includegraphics[width=1.1\textwidth]{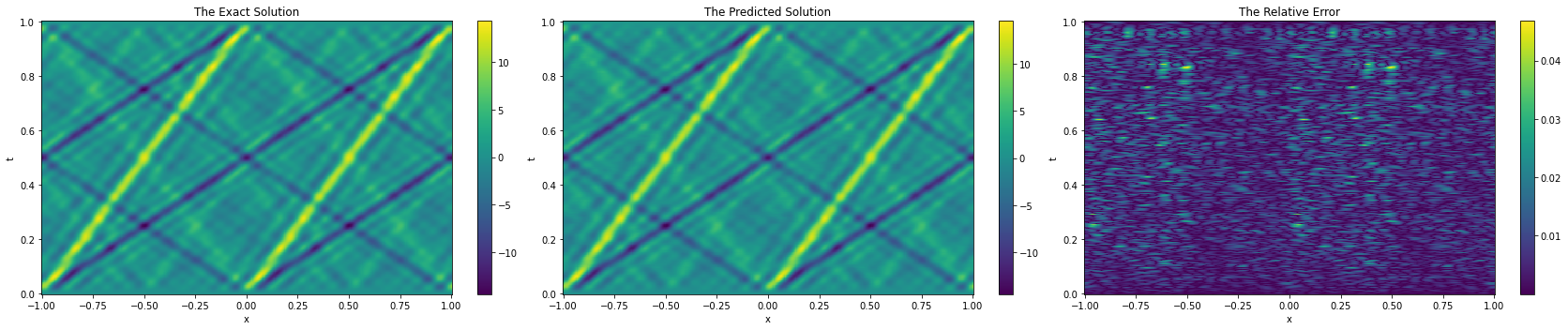}}
\caption{Random testing case after 500 epochs training for the simple case.
The maximum relative error is $5\%$.}%
\label{fig:1dmode10}%
\end{figure}

\noindent\textbf{Case 1: Constant wave-speed case} In this test, we assume the
wave speed is $2$, and the right hand side $f\left(  x \right) $ satisfies
\begin{equation}
f\left(  x,t \right)  = \sum_{i=0}^{2} c_{i}t ^{i} + \sum_{i=1}^{10}
2(c^{2}-1)i^{2}\pi^{2} \left(  a_{i}\cos(2 i \pi t)\cos\left(  2i \pi x
\right)  + b_{i}\sin(2 i \pi t) \sin(2 i \pi x) \right), \label{1drhs_test}%
\end{equation}
where $a_{i},b_{i},c_{i}$ are randomly sampled from $\left[  0,1 \right]  $.
The exact solutions are from integrating the formula based on Duhamel's
method. There are 21 basis based on $\ref{1drhs_test}$, thus, the input size
of the causality DeepONet should be $21 \times N_{t}$, where $N_{t}$ is the
number of time steps. We set $N_{t} = 400$ for all other cases in this paper.
We build the branch part and trunk part of the causality DeepONet as $4$-layer
neural network with $128$ hidden neurons each layer and $4$-layer neural
network with $100$ hidden neurons each layer respectively, both of whose
outputs dimension is $500$ and shrinking to dimension $N_{x} \times N_{t}$ by
inner product with the $x$-basis. The learning rates during training including
the following cases are $10^{-4}$. The batch size we choose for all the
testing case is 20. The results trained 500 epochs with 1000 cases are shown
in Figure \ref{fig:1dmode10}.  

\bigskip\noindent\textbf{Case 2: Variable wave-speed case}

\begin{figure}[ptbh]
\makebox[\textwidth][c]{\includegraphics[width=1.1\textwidth]{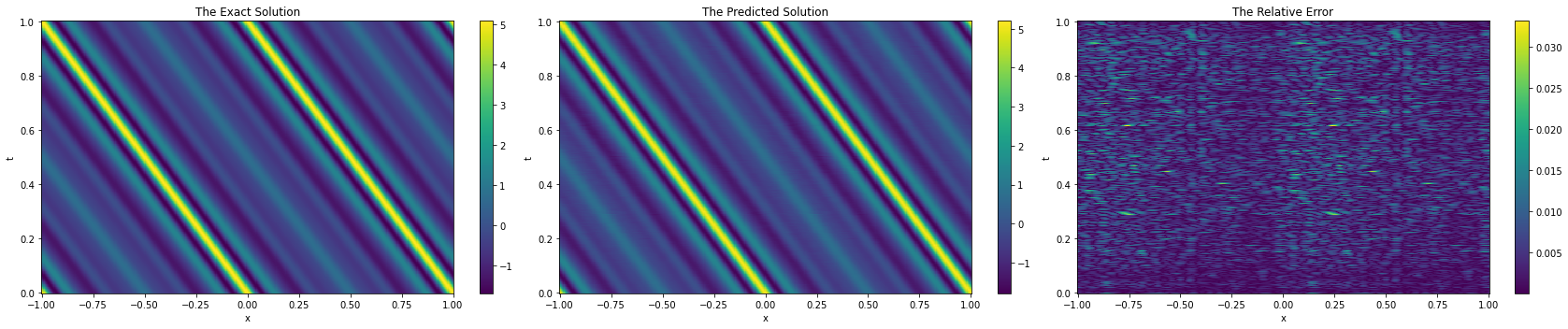}}
\caption{Random testing case after 500 epochs training for the variable wave
speed case. The maximum relative error is $4.2\pm2.3\%$.}%
\label{fig:1dcxtmode10}%
\end{figure}

\begin{figure}[ptbhptbh]
\makebox[\textwidth][c]{\includegraphics[width=1.1\textwidth]{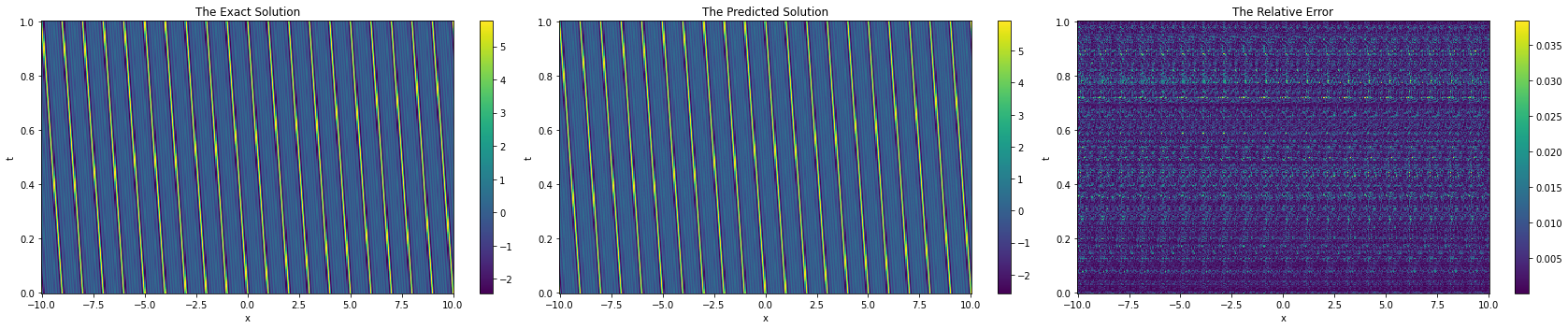}}
\caption{Random testing case after 500 epochs training for the variable wave
speed case with random large domain on $x \in\left[  -10,10 \right]  $. The
maximum relative error is $4.0 \pm1.4 \%$.}%
\label{fig:1dcxtmode10_large}%
\end{figure}

In this case we consider the case the wave speed is a function of $x,t$
\begin{equation}
c\left(  x,t \right)  = \sqrt{\cos\left(  2 n \pi\left(  t+x \right)  \right)
+ 1},\label{wavespeed}%
\end{equation}
and the exact solution is given by
\begin{equation}
u\left(  x,t \right)  = \sum_{m=1}^{M} c_{m} \cos\left(  2 m \pi\left(  t + x
\right),  \right) \label{1dvsol}%
\end{equation}
thus, the corresponding right hand side for (\ref{1dwave}) will be
\begin{equation}
\begin{aligned} f\left( x,t \right) = \sum_{m=1}^{M} c_{m} \left( 2m\pi \right)^2 &\left[ \cos \left( 2\pi \left( m+n \right)t \right) \cos\left( 2\pi\left( m+n \right)x \right) + \right. \\ & - \sin\left( 2\pi\left( m+n \right)t \right)\sin\left( 2\pi\left( m+n \right)x \right) \\ & + \cos\left( 2\pi\left( n-m \right)t \right)\cos\left( 2\pi\left( n-m \right)x \right)\\ & \left. - \sin\left( 2\pi \left( n-m \right)t \right)\sin\left( 2\pi\left( n-m \right)x \right) \right] . \end{aligned}\label{rhs_1dv}%
\end{equation}
The inputs for the CPOD-DeepONet are the coefficients of basis funciton of $x$
in $\left( \ref{rhs_1dv} \right) $, since there are $4M$ basis, the input size
of the causality DeepONet should be $4M \times N_{t}$. The causality DeepONet
is as the same shape as in Case 1 but with different input dimensions. In this
case we assume $M = 10$ and $n = 10$ with 400 equal spaced time steps, and
there are 400 random sampled points on x-direction for each time step. We have
similar results as shown in Figure \ref{fig:1dcxtmode10}. It should be mentioned that the trained CPOD-DeepONet
can in fact be used to predict solution outside the trained spatial domains as shown in Fig. \ref{fig:1dcxtmode10_large}.

 

\section{A Recursive DeepPropNet for Learning Evolution PDE Operators}
\label{sec:dp}
Following the semi-group formulation of evolution PDEs, the DeepPropNet computes block
by block along the time direction recursively. The initial conditions for the PDEs will be included
 together with the force functions for the time block as an input for DeepPropNet to assure the well-posedness of the problem. And, the DeepPropNet is then constructed recursively with a Causality DeepONet with POD to learn the PDE evolution operator
over long time interval.

To simplify notation, we
denote the initial conditions on the solution $u$ and its velocity $v$ as
\begin{equation}
u_{d}\left(  {x}\right)  =\left[  u_{0}\left(  {x}\right)
,v_{0}({x}),u_{1}\left(  {x}\right)  ,v_{1}({x}%
),\ldots,u_{n}\left(  {x}\right)  ,v_{n}({x})\right],
\end{equation}
and the
DeepPropNet $P_{\theta}$. Following the notations of Causality DeepONet with
POD in equation $\left(  \ref{cpod}\right)  $, we define the DeepPropNet as%
\begin{equation}
\mathcal{P}_{\theta}(\overrightarrow{f})(x,t) \sim\left\{  \overrightarrow
{\sigma}_{br,v}\left(  \overrightarrow{u}_{d} \right) \odot\overrightarrow
{\sigma}_{trk,v}\left(  t \right)  + \overrightarrow{\sigma}_{br, c}\left(
\overrightarrow{f} \right) \odot\overrightarrow{\sigma}_{trk, c}\left(  t
\right) \right\}  \cdot\overrightarrow{\sigma}_{basis}\left(  x \right),
\label{DP_CPOD}%
\end{equation}
where $\overrightarrow{\sigma}_{br,c}\left(  \overrightarrow{f} \right) $ is
the branch net $\left(  \ref{branch_cpod} \right) $, $\overrightarrow{\sigma
}_{trk,c}\left(  t \right) $ is the corresponding trunk net as shown in
equation $\left(  \ref{trunk_cpod} \right) $. Likewise, $\overrightarrow
{\sigma}_{br,v}\left(  \overrightarrow{u}_{d} \right)  \text{ and }
\overrightarrow{\sigma}_{trk,v}\left(  t \right) $ are the trunk net and
branch net of the vanilla DeepONet,
\begin{equation}
\begin{aligned} \overrightarrow{\sigma}_{br, v, i}\left( \overrightarrow{u}_{d} \right) = \sigma _{b}&\left( \sum_{n=0}^{N} \sum_{j=1}^{N}\xi_{i,j}^{n}u_n\left( s_{j}\right) + b_{i}\right) ,\end{aligned}\label{branch_c}%
\end{equation}
the corresponding $i$-th output of trunk net is
\begin{equation}
\overrightarrow{\sigma}_{trk,v,i}\left(  t \right)  = \sigma_{trk}\left(
\omega_{i}\cdot t+\zeta_{i}\right) .\label{trunk_c}%
\end{equation}
The notation is abused to some extent in above equations, but the $\sigma_{trk}$'s and
$\sigma_{b}$'s on the right hand side means different neural networks. The
definition of $\overrightarrow{\sigma}_{basis}$ is as the same as $\left(
\ref{basis} \right) $. Note the summation of $\overrightarrow{\sigma}%
_{br,v}\left(  \overrightarrow{u}_{d} \right) \odot\overrightarrow{\sigma
}_{trk,v}\left(  t \right) $ and $\overrightarrow{\sigma}_{br, c}\left(
\overrightarrow{f} \right) \odot\overrightarrow{\sigma}_{trk, c}\left(  t
\right) $ in the equation $\left(  \ref{DP_CPOD} \right) $ is inspired from
the superposition of linear waves from different sources and other format could be
explored for nonlinear problems.

\bigskip\noindent$\bullet$ \textbf{Recursive Formulation of DeepPropNet:} To
predict the wave field in the time block $\left[  t_{0},t_{1} \right]  $ at
beginning, the Deep Propagator DeepPropNet is to learn the mapping
\begin{equation}
P_{\theta}: \left[  u_{d}\left(  {x}_{j} \right) ,f\left(
{x}_{\xi},s \right) , t_{0} \le s \le t_{1} \right]  \mapsto
u\left(  {x},t \right) \quad t_{0} \le t \le t_{1}, {x}
\in\mathbb{R}^{n}.\label{DP_1}%
\end{equation}

Once the Deep Propagator is learned, the wave field in the next time block
$\left[  t_{1},t_{2} \right]  $ could be predicted by the Deep Propagator by
using the initial propagator as follows
\begin{equation}
u\left(  x,t \right)  = P_{\theta}\left[  P_{\theta}\left(  {x}_{j}
,t_{1}\right) ,{\dot P}_{\theta}\left(  {x}_{j}
,t_{1}\right) ,f\left(  {x}_{\xi},s \right) , t_{1} \le s \le t_{2}
\right]  \mapsto u\left(  {x},t \right) \quad t_{1} \le t \le
t_{2}, {x} \in\mathbb{R}^{n},\label{DP_2}%
\end{equation}
where the initial condition is replaced by the prediction of DeepPropNet at
time $t_{1}$.

The resulting propagator for the time block $[t_{1}, t_{2}]$ then will be
again to be used to provide the initial condition for $t=t_{2}$ and the same
initial propagator network for the time period $[t_{0}, t_{1}]$. This
procedure can be carried on recursively until the whole time period $[0, T]$
is covered,a global propagator network is thus obtained. Since it is like the
initial propagator solver tracks the waves and propagates with the solutions
along time direction and the propagator itself is a deep neural network,
we call it Deep Propagator.

The schematics of the Deep Propagator DeepPropNet are shown in Algorithm \ref{alg:DP}.
This is another crux to solve high dimensional problem. As illustrated in
section \ref{sec:evoop}, the input size could explode since the global dependence
of source term, by solving it block by block recursively combining with the
memory-efficient trick from separation of variables, learning high dimensional evolutioin operators is
manageable now.

\textbf{Loss function} Similar to the loss function in section \ref{sec:cpod}, we will define the loss function of DeepPropNet block by block in time. Given a batch size $\mathcal{N}$ for training process, number of blocks \(N_{b}\) and
the total number of the test records $N$, the loss function is defined as
\begin{equation}
\mathcal{L}oss\left(  \theta\right)  =\frac{1}{\mathcal{N}}\sum_{i=1}%
^{\mathcal{N}} \sum_{n=1}^{N_{b}} \sum_{j=N_{t,n}}^{N_{t,n+1}}\sum_{k=1}^{N_{{x}}} \left(
\mathcal{G}\left(  \overrightarrow{f} \right) \left(  {x}_{k}%
,t_{j}\right)  -y_{ijk}\right)  ^{2}, \label{loss}%
\end{equation}
where the \(n\)-th block starts from \(t_{N_{t,n}}\) but ends in \(t_{N_{t,n+1}}\) and $N_{x}$ is the number of points
on ${x}$ direction.

Likewise, the mean of the relative L2 error is
considered as the relatvie L2 error in the total time-spatial domain
\begin{equation}
\mathcal{L}oss_{2}\left(  \theta\right)  = \frac{1}{\mathcal{N}}\sum
_{i=1}^{\mathcal{N}}\frac{\sum_{i=1}^{N_{b}} \left\|  \mathcal{G}\left(  \overrightarrow{f_{i}}
\right) \left(  {x},t \right)  - u_{i}\left(  {x},t
\right)  \right\|_{2,\Omega_{i}}}{\sum_{i=1}^{N_{b}} a_i z^i\left\|  u_{i}\left(  {x},t \right)
\right\|_{2,\Omega_{i}} },
\end{equation}
where $\Omega_i$ is the i-th block of time and spatial domain.

\begin{algorithm}
\caption{The schematics of Deep Propagator}\label{alg:DP}
\SetKwInOut{Input}{Input}
\SetKwInOut{Output}{Output}
\Input{\(
\left\{\begin{aligned}
&t \text{: Time step}\\
&{x} \text{: Spatial Coordinates}\\
&u_{d}\left( {x} \right) \text{: Initial conditions}\\
&\left[ a_0\left( s \right), a_1\left( s \right), b_1\left( s \right),\ldots ,a_{N}\left( s \right),b_{N}\left( s \right) \right] \text{: Coefficients of right hand side} \\
&N \text{: Number of blocks at time domain}
.\end{aligned}\right.
\)   }
\Output{$u\left( {x},t \right)$}
\(u_0\left( {x} \right) \gets u_{d}\left( {x} \right)\)
\For{\(i = 0; i<N; i++\)}{
\(u\left( {x},t \right) \gets \text{DeepPropaNet}\left( t, u_0\left( {x}_{j} \right), \left[ a_0\left( s \right), a_1\left( s \right), b_1\left( s \right),\ldots ,a_{N}\left( s \right),b_{N}\left( s \right) \right] \right)\)\;
$u_0\left( {x} \right) \gets u\left( {x}, t_{i+1} \right) $\;
}
\end{algorithm}
\begin{algorithm}
\caption{The algorithm of DeepPropaNet}\label{alg:DeepPropaNet}
\SetKwInOut{Input}{Input}
\SetKwInOut{Output}{Output}
\Input{\(
\left\{\begin{aligned}
&t \text{: Time step}\\
&{x} \text{: Spatial Coordinates}\\
&u_{d}\left( {x} \right) \text{: Initial conditions}\\
&\left[ a_0\left( s \right), a_1\left( s \right), b_1\left( s \right),\ldots ,a_{N}\left( s \right),b_{N}\left( s \right) \right] \text{: Coefficients of right hand side}
\end{aligned}\right.
\)   }
\Output{$u\left( {x},t \right)$}
\(y_t \gets \text{CDeepONet}\left( t, \left[ a_0\left( s \right), a_1\left( s \right), b_1\left( s \right),\ldots ,a_{N}\left( s \right),b_{N}\left( s \right) \right] \right)\) \Comment*[r]{The output \(y_t\) has size of \([\text{batchsize} \times  \text{Nt} \times \text{NBx}]\)}
\(y_i \gets \text{DeepONet}\left( t, u_{d}\left( {x}_{j} \right) \right)\) \Comment*[r]{DeepONet is the vanilla DeepONet. The output \(y_i\) has size of \([\text{batchsize} \times  \text{Nt} \times \text{NBx}]\)}
$y_x \gets \left[ 1, \cos\left( 2\pi {x} \right), \sin\left( 2\pi {x} \right), \ldots \cos\left( 2 M \pi {x} \right), \sin\left( 2 M {x} \right)\right] $\Comment*[r]{The output \(y_x\) has size of \([\text{batchsize} \times  \text{Nx} \times \text{NBx}]\)}
$u\left( x,t \right) \gets \text{Einsum}\left( \text{'bxn,btn} \to \text{ bxt'}, y_{x}, y_{t}+y_{i} \right) $
\end{algorithm}

\section{Numerical Results of DeepPropNet}
\label{sec:results}
\begin{figure}[ptbh]
\makebox[\textwidth][c]{\includegraphics[width=1.1\textwidth]{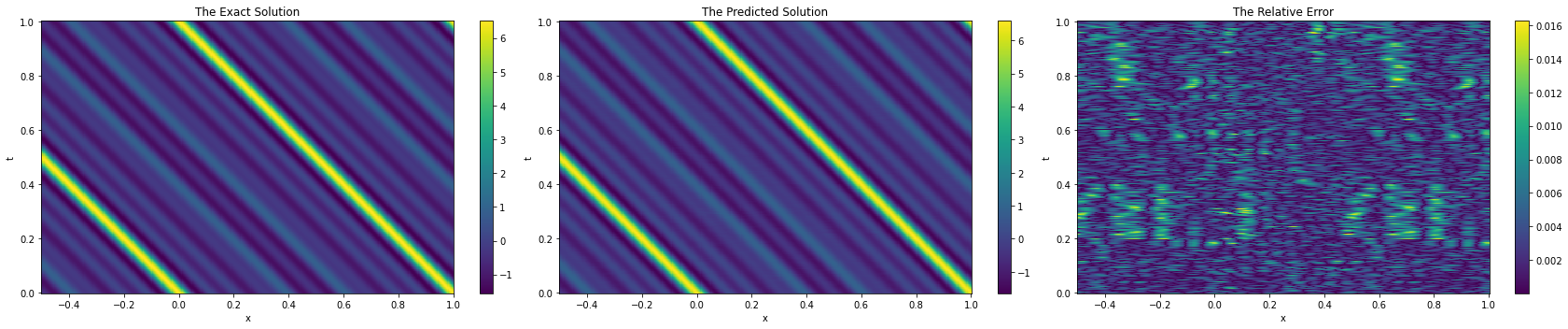}}
\caption{One of the random testing cases after 500 epochs training for the variable wave
speed case with DeepPropNet whose initial condition is given by true data
in $x \in\left[  -0.5,1 \right]  $. The maximum relative error in the whole
domain is $3.4\pm2.0 \%$. (Left) Exact, (Middle) Prediction, (Right) Relative Error.}%
\label{fig:dp_1dcxtmode10_true}%
\end{figure}

The equation and right hand sides are as the same as the non-constant wave
speed case $\left(  \ref{wavespeed} \right) -\left(  \ref{rhs_1dv} \right) $,
but we split 5 equal sized blocks in time domain and for each block, there are 80 sampled time
locations for training. The causality DeepONet is as the same shape as in Case 2 but with
different input dimensions since the dependence of time steps shrinks. The vanilla DeepONet
is a composition of a $4$-layer fully connected neural network with $128$
hidden neurons each layer as branch net and a $4$-layer fully connected neural
network with $100$ hidden neurons each layer as trunk net, both of whose
output dimensions are also $500$. The initial data for each block will be given by true data. Figure \ref{fig:dp_1dcxtmode10_true} shows the prediction
of the exact solution up to time $t=1$ with a maximum relative error of $3.4\pm2.0 \%$.



\begin{figure}[ptbh]
\makebox[\textwidth][c]{\includegraphics[width=1.1\textwidth]{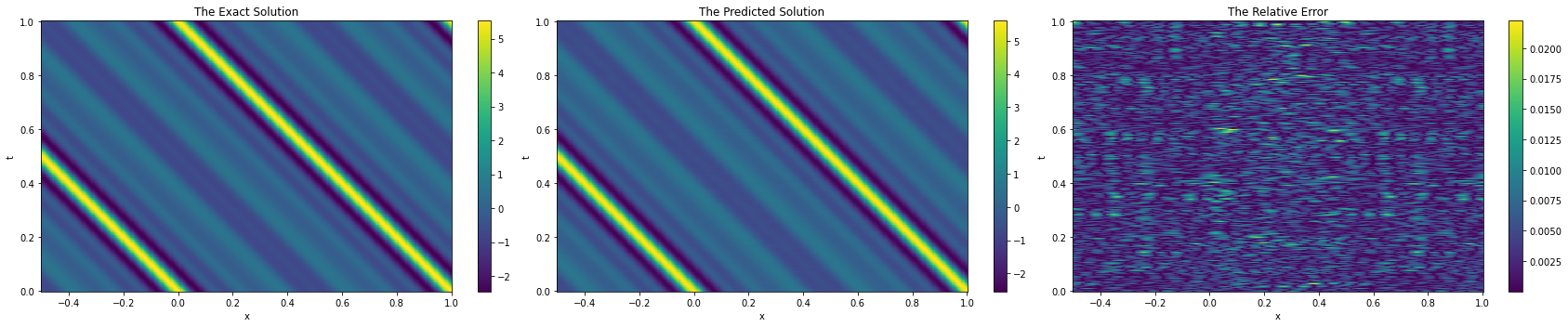}}
\caption{One of the random testing cases after 500 epochs training for the variable wave
speed case with DeepPropNet whose initial condition is given by the
prediction in $x \in\left[  -0.5,1 \right]  $. The maximum relative error in
the whole domain is $3.2\pm1.9 \%$. (Left) Exact, (Middle) Prediction, (Right) Relative Error.}%
\label{fig:dp_1dcxtmode10_pred}%
\end{figure}
\begin{figure}[ptb]
\centering  \begin{subfigure}[b]{0.45\textwidth}
\centering
\includegraphics[width=\textwidth]{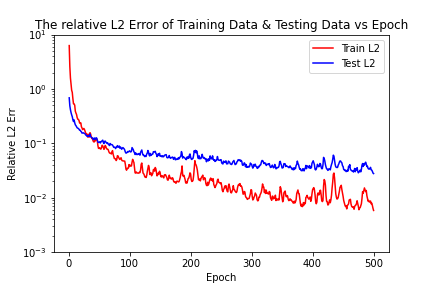}
\caption{The case using prediction \(P_{\theta}\left[ u_{0}\left( x_j \right), v_{0}\left( x_j \right), f \right]\left( x, t_1 \right)\) as initial condition for each time block.}
\label{fig:dp_lossPred}
\end{subfigure}
\hfill\begin{subfigure}[b]{0.45\textwidth}
\centering
\includegraphics[width=\textwidth]{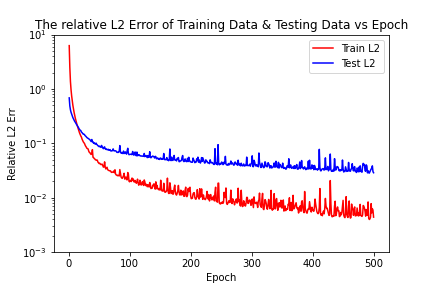}
\caption{The case using exact solution \(u\left( {x},t_{i} \right)\) for initial condition for each time block.}
\label{fig:dp_lossExact}
\end{subfigure}
\caption{Relative L2 error evolution during training}%
\label{fig:dp_loss}%
\end{figure}

\medskip\noindent\textbf{Training with DeepPropNet Prediction as Initial Condition at $t_i$.} There
are two strategy to provide initial conditions for each block, either provide
the exact solution $u\left(  x,t_{i} \right) $ or the prediction $P_{\theta
}\left[  u_{0}\left(  x_{j} \right) , v_{0}\left(  x_{j} \right), f \right] \left(  x, t_{1} \right) $. We
compare the relative L2 error evolution of these $2$ different process and as
shown in Figure \ref{fig:dp_loss}, they have similar loss convergence. The
results are shown in Figure \ref{fig:dp_1dcxtmode10_pred} and Figure
\ref{fig:dp_1dcxtmode10_true} respectively. This is a natural test about the
capability of the interpolation in training domain of the neural network, but
it offers a convenience that there is no further need to prepare exact data on
grids for initial conditions of each time blocks.

\section{Conclusion and Future works}
\label{sec:future}
In this paper, we proposed the recursive way to construct DeepPropNet - a DNN
propagator for evolution system over large time period by using a single
building block propagator over a small time period, thus reducing the overall
complexity and size of the neural network required. For the design of the DeepPropNet, we
also extended the Causality DeepONet with POD with specific basis
to alleviate the memory burden for large spatial variables. By seperately handling the spacial and temporal domain, we gain not only the memory efficiency but also the training boost.
The preliminary numerical results have shown the feasibility of this recursive DeepPropNet in predicting the time evolution of wave propagations.

The proposed DeepPropNet here is based on a supervised learning approach where
the data can be generated by a separate numerical methods or observation data
or analytical solution when available. In theory, we could also use a
unsupervised learning procedure to train the DeepPropNet by using the
following residue of the PDEs as the loss function.
\begin{equation}
\mathcal{L}\text{oss} = \int_{\Omega} \mid\mid\mathcal{L}_{t} P_{\theta
}\left[ u_{0}, v_0,  f \right]  - \mathcal{L}_{x}P_{\theta}\left[ u_{0}, v_0, f \right]
- f\mid\mid^{2} dxdt.\label{uloss}%
\end{equation}

It is also natural to extend this framework to learn high frequency problem, since the spatial oscillating terms will be handled explicitly by the basis, but the temporal oscillating terms could be handled by the Causality DeepONet, based on our experience in the work \cite{cdeepOnet}.

Future work will be conducted on more complex evolution systems, the
unsupervised training, training over partial time domain, highly oscillating problems and higher dimensional problems as well as initial boundary value
problems.

\end{document}